# HITTING PROBABILITIES IN A MARKOV ADDITIVE PROCESS WITH LINEAR MOVEMENTS AND UPWARD JUMPS: APPLICATIONS TO RISK AND QUEUEING PROCESSES

By Masakiyo Miyazawa[1]

*Tokyo University of Science*

Motivated by a risk process with positive and negative premium rates, we consider a real-valued Markov additive process with finitely many background states. This additive process linearly increases or decreases while the background state is unchanged, and may have upward jumps at the transition instants of the background state. It is known that the hitting probabilities of this additive process at lower levels have a matrix exponential form. We here study the hitting probabilities at upper levels, which do not have a matrix exponential form in general. These probabilities give the ruin probabilities in the terminology of the risk process. Our major interests are in their analytic expressions and their asymptotic behavior when the hitting level goes to infinity under light tail conditions on the jump sizes. To derive those results, we use a certain duality on the hitting probabilities, which may have an independent interest because it does not need any Markovian assumption.

**1. Introduction.** Consider a risk process in which premium rates, claim arrivals and claim sizes depend on a background state, which is governed by a continuous time Markov chain. The premium rate is usually assumed to be positive, and can be reduced to unit by a random time change. Suppose that the premium includes returns from investments. Then, it may happen that the premium rate takes negative values since the returns may be negative. This motivates us to consider a risk process with positive and negative premium rates. Our primary interest is in an asymptotic behavior of the ruin probability when the initial reserve goes to infinity, provided the claim size distributions have light tails, that is, decrease exponentially fast.

Received May 2002; revised April 2003.
[1]Supported in part by JSPS Grant 13680532.
*AMS 2000 subject classifications.* Primary 90B22, 60K25; secondary 60K20, 60G55.
*Key words and phrases.* Risk process, Markov additive process, hitting probability, decay rate, Markov renewal theorem, stationary marked point process, duality.

Similarly to Asmussen (2000), we model the risk process as a Markov additive process $Y(t)$ with upward jumps and $Y(0) = 0$. Then, for $x \geq 0$, $x - Y(t)$ is a risk process with initial reserve $x$, and ruin occurs when $Y(t)$ hits level $x$. We assume the following setting. The background continuous-time Markov chain has a finite state space. The additive component linearly increases or decreases while the background state is unchanged. At the transition instants of the background Markov chain, the additive component may have upward jumps whose sizes may depend on the background states before and after the transitions. This additive process was recently studied by Takada (2001) and Miyazawa and Takada (2002) for extending the Markov modulated fluid queue [see, e.g., Asmussen (1995) and Rogers (1994)].

A key observation in Miyazawa and Takada (2002) and Takada (2001) is that the hitting probabilities in the downward direction have a matrix exponential form, which follows from the fact that the additive process is skip free in this direction. The result generalizes the matrix exponential forms in the literature [see, e.g., Asmussen (1995, 1995), Rogers (1994) and Takine (2001)]. However, we can not directly apply this result to the ruin probability since this requires the hitting probabilities in the upward direction.

In this paper we consider these hitting probabilities. We represent the hitting probabilities as a Markov renewal function. Then, their asymptotic behaviors are studied through the Markov renewal theorem. This renewal approach for the asymptotic behavior is rather classical, but has not been so popular for the Markov additive process. A typical approach is the change of measures based on martingales [see, e.g., Asmussen (2000)]. The large deviation technique has also been used. The latter can be applied for a larger class of additive processes, but is less informative. In this paper we show that the Markov renewal approach is yet powerful, in particular, it provides the prefactor of the decay function.

To fully utilize this approach, we need to get the Markov renewal kernel in a closed form, which is the conditional joint distribution of the ladder height and associated background state given an initial background state. To this end, we use a dual additive process similarly to Asmussen (1991) and Sengupta (1989). This dual process is obtained from the Markov additive process by time reversal and changing the sign of the additive component. However, our approach is different from those in Asmussen (1991) and Sengupta (1989), which use the occupation measure. We directly derive the joint distribution, including the ladder epoch, in terms of the dual process, which does not require any drift condition on the additive process. Furthermore, we do not need any Markovian assumption, but only need stationarity. This result may be of independent interests, so we derive it under the stationary regime. We then specialize it to the Markov additive process and perform detailed computations.

This paper is made up by six sections. In Section 2 we introduce the additive process with upward jumps under the stationary regime, and compute the ascending ladder height distribution in terms of the dual process. In Section 3 we specialize this result to the Markov additive processes. In Section 4 we derive a Markov renewal equation for the hitting probabilities. In Section 5 those probabilities are shown to have asymptotically exponential decay. In Section 6 we discuss applications to risk processes and fluid queues.

**2. The additive process under the stationary framework.** In this section we introduce a stationary additive process and consider its ascending ladder height, not using any Markovian assumption. Let $\{(t_n, M_n, A_n); n \in \mathbb{Z}\}$ be a stationary marked point process, such that the counting measure $N$ of $\{t_n\}$ has a finite intensity $\lambda$ and the mark $(M_n, A_n)$ takes values in $S \times \mathbb{R}_+$, where $S$ is a countable set, $\mathbb{Z}$ is the set of all integers and $\mathbb{R}_+ = [0, \infty)$. In our application $S$ is finite, but this does not matter in this section. In the terminology of a risk process, $t_n$ is the arrival time of a claim if $A_n$ is positive and $M_n$ is a background state. Note that not all $t_n$ represent claim arrivals since $A_n = 0$ is allowed. This gives a great flexibility when we assume $M_n$ to be Markov. Namely, the claim arrivals are subject to the Markovian arrival process due to Neuts (1989) [see Asmussen (2000) for applications to a risk process].

Let $M(t)$ and $A(t)$ be continuous-time processes such that $M(t) = M_n$ and $A(t) = A_n$, respectively, for $t \in [t_n, t_{n+1})$. As usual, we can choose a canonical probability space $(\Omega, \mathcal{F}, P)$ with a shift operator group $\{\theta_t\}$ on $\Omega$ such that $P(\{\theta_t(\omega) \in D\}) = P(D)$ for all $D \in \mathcal{F}$, and, for $\omega \in \Omega, s \in \mathbb{R} \equiv (-\infty, +\infty)$,

$$N(B)(\theta_s(\omega)) = N(B + s)(\omega), \qquad B \in \mathcal{B}(\mathbb{R}),$$

$$M(t)(\theta_s(\omega)) = M(t + s)(\omega), \quad A(t)(\theta_s(\omega)) = A(t + s)(\omega), \qquad t \in \mathbb{R},$$

where $B + s = \{u + s; u \in B\}$. We shall use abbreviated notation such as $M(t) \circ \theta_s$ for $M(t)(\theta_s(\omega))$.

We next define an additive process. Let $v(i)$ be a function from $S$ to $\mathbb{R} \setminus \{0\}$. For simplicity, we exclude the case $v(i) = 0$, but this is not essential. Then, $Y(t)$ is defined as

$$(2.1) \qquad Y(t) = \begin{cases} \int_{0+}^{t+} v(M(u)) \, du + \int_{0+}^{t+} A(u) N(du), & t \geq 0, \\ -\int_{t-}^{0+} v(M(u)) \, du - \int_{t-}^{0+} A(u) N(du), & t < 0. \end{cases}$$

Thus $Y(t)$ is the additive process that linearly changes while $M(t)$ is constant and may have jumps when it changes. Note that $Y(t)$, $M(t)$ and $A(t)$ are right continuous and $Y(0) = 0$. By definition, $Y(t)$ has stationary increments and satisfies $(Y(t) - Y(s)) \circ \theta_u = Y(t + u) - Y(s + u)$. In particular, $Y(0) = 0$ implies

$$Y(t) \circ \theta_u = Y(t + u) - Y(u). \tag{2.2}$$

We refer to $Y(t)$ as a stationary additive process with background process $M(t)$.

Let $\tau_x^+$ be the hitting time when $Y(t)$ gets into $[x, \infty)$ for $x \geq 0$. Namely,

$$\tau_x^+ = \inf\{t > 0; Y(t) \geq x\}.$$

Then, $Y(\tau_0^+)$ is called the ascending ladder height, which is zero when $Y(t)$ continuously crosses the origin. Note that $\tau_0^+$ and, hence, $Y(\tau_0^+)$ and $M(\tau_0^+)$ are defective (proper) if

$$E(Y(1)) < (>)0 \qquad \text{respectively.}$$

As is well known, this can be proved using Loynes' (1962) arguments. If $M(t)$ is a Markov process, $E(Y(1)) \geq 0$ implies that $\tau_0^+$ is proper [see Lemma 2.1 of Miyazawa and Takada (2002)]. In what follows, for $i \in S$, the event $\{M(\tau_0^+) = i\}$ means the event $\{\tau_0^+ < \infty, M(\tau_0^+) = i\}$.

The purpose of this section is to express the joint distribution at the ascending ladder epoch, that is, for $i, j \in S, u, x, z \geq 0$,

$$P(\tau_0^+ \leq u; M(0) = i, M(\tau_0^+ -) = j, M(\tau_0^+) = k, -Y(\tau_0^+) \leq z, Y(\tau_0^+) \leq x),$$

in terms of dual processes $\tilde{M}(t)$, $\tilde{A}(t)$ and $\tilde{Y}(t)$ defined below. In our applications we only need the case that $u = z = \infty$, but the joint distribution may be of independent interests, as we shall see.

Define dual processes $\tilde{M}(t)$ and $\tilde{Y}(t)$ by

$$\tilde{M}(t) = M(-t), \qquad \tilde{A}(t) = A(-t), \qquad \tilde{Y}(t) = -Y(-t) + Y(0-), \qquad t \in \mathbb{R}.$$

Note that $\tilde{Y}(t)$ is also a stationary additive process with upward jumps satisfying $\tilde{Y}(0+) = 0$. The term $Y(0-)$ is irrelevant in the definition of $\tilde{Y}(t)$ under the stationary setting, that is, with respect to $P$, since $P(Y(0) = Y(0-) = 0) = 1$. However, this is not the case under the conditional probability measure, given that there is a point of $N$ at the origin. Note that the dual process $\tilde{Y}(t)$ is not right continuous but left continuous. Thus, $\tilde{Y}(0) = \tilde{Y}(0-)$ may be negative. Similarly to $\tilde{Y}(t)$, $\tilde{M}(t)$ and $\tilde{A}(t)$ are left continuous. Furthermore, if $\tilde{M}(t) = i$, then $\tilde{Y}(t)$ changes at rate $v(i)$, that is, $Y(t)$ and $\tilde{Y}(t)$ have the same rate when their corresponding background states are identical. Let

$$\tilde{\tau}_w^- = \inf\{t > 0; \tilde{Y}(t) \leq w\}, \qquad w \leq 0.$$

That is, $\tilde{\tau}_w^-$ is the hitting time of $\tilde{Y}(t)$ at level $w \leq 0$. Note that $\tilde{Y}(t)$ is skip free in the downward direction. So the distributions of $\tilde{\tau}_w^-$ and $\tilde{M}(\tilde{\tau}_w^-)$ are easier to get when the background process $M(t)$ is Markov. This fact will be used.

Since $\lambda \equiv E(N(0,1]) < \infty$, we can define the Palm distribution with respect to the point process $N$, which is denoted by $P_N$. Let

$$\sigma^-(t) = \sup\{u < t; N([u,t)) > 0\},$$

$$\sigma^+(t) = \inf\{u > t; N([t,u)) > 0\}.$$

We also use their duals: $\tilde{\sigma}^-(t)$ and $\tilde{\sigma}^+(t)$, by changing $N$ to its dual $\tilde{N}$, where $\tilde{N}(B) = N(-B)$ for $B \in \mathcal{B}(\mathbb{R})$ in which $-B = \{x \in \mathbb{R}; -x \in B\}$. For example, $\tilde{\sigma}^-(t)$ is the last transition instant of $\tilde{M}$ before time $t$. Note that the dual processes under $P$ coincide with the original processes under $\tilde{P}$, which is a canonical probability measure for $(\tilde{N}(\cdot), \tilde{M}(t), \tilde{A}(t))$. For example, $E_N(\varphi(\tilde{\tau}_{-x}^-)) = \tilde{E}_{\tilde{N}}(\varphi(\tau_{-x}^-))$. Keeping these facts in mind, we present the following results, using the following notation:

$$S^- = \{i \in S; v(i) < 0\}, \qquad S^+ = \{i \in S; v(i) > 0\},$$

$$\mathbf{M}^-(t) = (M(\sigma^-(t)-), M(t)), \qquad \tilde{\mathbf{M}}^+(t) = (\tilde{M}(t), \tilde{M}(\tilde{\sigma}^+(t)+)).$$

Note that $\mathbf{M}^-(t) = (M(t-), M(t))$ if $N(\{t\}) = 1$.

LEMMA 2.1. *For a nonnegative measurable $\varphi$, $i \in S^-$, $j,k \in S$ and $x,z \geq 0$,*

(2.3)
$$-v(i)E(\varphi(\tau_0^+); M(0) = i, \mathbf{M}^-(\tau_0^+) = (j,k),$$
$$-Y(\tau_0^+-) \leq z, 0 < Y(\tau_0^+) \leq x)$$
$$= \lambda \int_0^z E_N(\varphi(\tilde{\tau}_{-w}^-); \tilde{\mathbf{M}}^+(0) = (k,j), \tilde{M}(\tilde{\tau}_{-w}^-) = i,$$
$$0 < \tilde{A}(0) - w \leq x) \, dw,$$

*and, using the notation $\tilde{T} = \tilde{Y}(0) - \tilde{Y}(\tilde{\sigma}^-(0)+)$, for $i \in S^-$, $k \in S^+$,*

(2.4)
$$-v(i)E(\varphi(\tau_0^+); M(0) = i, \mathbf{M}^-(\tau_0^+) = (j,k), Y(\tau_0^+) = 0)$$
$$= \lambda \int_0^\infty E_N(\varphi(\tilde{\tau}_{-w}^-); \tilde{\mathbf{M}}^+(0) = (k,j), \tilde{M}(\tilde{\tau}_{-w}^-) = i,$$
$$\tilde{A}(0) < w < \tilde{T} + y) \, dw,$$

*where $E(X; D) = E(X \mathbb{1}_D)$ for a random variable $X$ and an event $D$.*

PROOF. We first consider the case that $Y(\tau_0^+) > 0$. In this case the additive process $Y(t)$ only up crosses level 0 by a jump. Hence, we have

$$E(\varphi(\tau_0^+); M(0) = i, \mathbf{M}^-(\tau_0^+) = (j,k), -Y(\tau_0^+) < z, 0 < Y(\tau_0^+) \leq x)$$

$$\begin{aligned}
&= E\bigg(\int_0^\infty \varphi(t)\mathbb{1}(M(0)=i, \mathbf{M}^-(t)=(j,k),\\
&\qquad\qquad Y(u)<0, u\in(0,t),\\
&\qquad\qquad -Y(t-)\le z, Y(t)\in(0,x])N(dt)\bigg),
\end{aligned}$$
(2.5)

where $\mathbb{1}(\cdot)$ is the indicator function of the statement "·". To (2.5), we apply Campbell's formula [see, e.g., (3.3.1) of Baccelli and Brémaud (2002)], which is given for a random function $g(t)$ by

$$E\bigg(\int_0^\infty g(t)\circ\theta_t N(dt)\bigg) = \lambda E_N\bigg(\int_0^\infty g(t)\,dt\bigg).$$

To this end, we apply $\theta_{-t}$ to the events in the indicator function of (2.5). By (2.2),

$$Y(t)\circ\theta_{-t} = -Y(-t) = \tilde{Y}(t) - \tilde{Y}(0).$$

Thus, using the fact that $Y(0+)=0$, we have

$$\begin{aligned}
\{M(0)=i, M(t-)&=j, M(t)=k\}\circ\theta_{-t}\\
&= \{\tilde{M}(0)=k, \tilde{M}(0+)=j, \tilde{M}(t)=i\},\\
\{Y(u)<0, u\in(0,t)\}&\circ\theta_{-t}\\
&= \{\tilde{Y}(t)\le \tilde{Y}(u), u\in(0,t)\},\\
\{-Y(t-)\le z, Y(t)&\in(0,x]\}\circ\theta_{-t}\\
&= \{-\tilde{Y}(t)\le z, -\tilde{Y}(0)-x\le -\tilde{Y}(t)\le -\tilde{Y}(0)\}.
\end{aligned}$$

Hence, applying Campbell's formula to the right-hand side of (2.5) yields

$$\begin{aligned}
E(\varphi(\tau_0^+); M(0)&=i, \mathbf{M}^-(\tau_0^+)=(j,k), 0<Y(\tau_0^+)\le x)\\
&= \lambda E_N\bigg(\int_0^\infty \varphi(\tilde{\tau}^-_{\tilde{Y}(t)})\mathbb{1}(\tilde{\mathbf{M}}^+(0)=(k,j), \tilde{M}(\tilde{\tau}^-_{\tilde{Y}(t)})=i,\\
&\qquad\qquad \tilde{Y}(t)<\tilde{Y}(u), u\in(0,t), -\tilde{Y}(t)\le z,\\
&\qquad\qquad -\tilde{Y}(0)-x\le -\tilde{Y}(t)<-\tilde{Y}(0))\,dt\bigg),
\end{aligned}$$
(2.6)

where the fact that $t=\tilde{\tau}^-_{\tilde{Y}(t)}$ when $t$ is the descending ladder epoch is used. We now change variable $t$ in the integral in (2.6). Let

$$w(t) = -\inf_{0<u<t}\tilde{Y}(u).$$

Clearly, $w(t)$ is nondecreasing, and $w(t) = -\tilde{Y}(t)$ when $w(t)$ increases because the dual process $\tilde{Y}(t)$ is left continuous and skip free in the downward direction. So, we have

$$w'(t) = -v(\tilde{M}(t))\mathbb{1}(\tilde{Y}(t) < \tilde{Y}(u), u \in (0,t)).$$

Hence, changing variables from $t$ to $w = w(t)$ in (2.6) and using the fact that $\tilde{A}(0) = -\tilde{Y}(0)$, we obtain (2.3). We next consider the case that $k \in S^+$ and $Y(\tau_0^+) = 0$. Similarly to (2.5), we have

(2.7)
$$E(\varphi(\tau_0^+); M(0) = i, \mathbf{M}^-(\tau_0^+) = (j,k), Y(\tau_0^+) = 0)$$
$$= E\bigg(\int_0^\infty \varphi(t)\mathbb{1}(M(0) = i, \mathbf{M}^-(t) = (j,k),$$
$$Y(u) < 0, u \in (0,t),$$
$$Y(t) \le 0 < Y(\sigma^+(t)-))N(dt)\bigg).$$

Hence, we can apply similar arguments as for the case $Y(\tau_0^+) > 0$, but we need to replace the event $\{0 < Y(t) \le x\} \circ \theta_{-t}$ by

$$\{Y(t) \le 0 < Y(\sigma^+(t)-)\} \circ \theta_{-t}$$
$$= \{-Y(-t) \le 0 < Y(\sigma^+(0)-) - Y(-t)\}$$
$$= \{\tilde{A}(0) \le w < \tilde{Y}(0) - \tilde{Y}(\tilde{\sigma}^-(0)+) + \tilde{A}(0)\},$$

for $w = -\tilde{Y}(t)$. Hence, (2.4) follows. □

REMARK 2.1. (a) A key fact in Lemma 2.1 is that the ladder height distribution, given the initial state, can be obtained via the hitting times in the skip free direction.

(b) Lemma 2.1 does not need any drift condition on $Y(t)$. So, the distribution in the left-hand side of (2.3), with $\varphi \equiv 1$, may or may not be defective.

(c) Consider the special case that $N$ is Poisson and $v(1) = -1$ with $S = \{1\}$. This constitutes the classical risk process. Then, (2.3) is compatible with standard results, such as Theorem 2.2 in Chapter III of Asmussen (2000), provided $E(Y(1)) < 0$.

In the rest of this section we shall give further remarks and results on Lemma 2.1. We first note the interesting feature of (2.3), that the sample path of the forward process is traced back under the Palm distribution. Figure 1 illustrates this sample path behavior. This fact is also related to the following conditional trace back for the conventional Markov modulated risk process, that is, in the Markovian setting with $v \equiv -1$.

$$
\begin{aligned}
(2.8) \quad & P(\tau_0^+ \leq t | M(0) = i, \mathbf{M}^-(\tau_0^+) = (j,k), -Y(\tau_0^+-) = z) \\
& = P(\tilde{\tau}_{-z}^- \leq t | \tilde{\mathbf{M}}^+(0) = (k,j), \tilde{M}(\tilde{\tau}_{-z}^-) = i), \qquad t, z \geq 0.
\end{aligned}
$$

See, for example, (2.1) of Asmussen and Højgaard (1996) and Proposition 2.2 of Asmussen and Klüppelberg (1996) for the derivation of (2.8). Let us show how (2.8) is obtained in our approach. From (2.3), we have, for $i \in S^-$, $j, k \in S$ and $t, z > 0$,

$$
\begin{aligned}
(2.9) \quad & -v(i) P(\tau_0^+ \leq t, M(0) = i, \mathbf{M}^-(\tau_0^+) = (j,k), 0 < -Y(\tau_0^+-) \leq z) \\
& = \lambda \int_0^z P_N(\tilde{\tau}_{-w}^- \leq t, \tilde{\mathbf{M}}^+(0) = (k,j), \tilde{M}(\tilde{\tau}_{-w}^-) = i, w < \tilde{A}(0)) \, dw.
\end{aligned}
$$

Hence, taking the derivative of (2.9) with respect to $z$ and dividing by this derivative with $t = \infty$, we get, for $t \geq 0$,

$$
\begin{aligned}
(2.10) \quad & P(\tau_0^+ \leq t | M(0) = i, \mathbf{M}^-(\tau_0^+) = (j,k), -Y(\tau_0^+-) = z) \\
& = P_N(\tilde{\tau}_{-z}^- \leq t | \tilde{\mathbf{M}}^+(0) = (k,j), \tilde{M}(\tilde{\tau}_{-z}^-) = i, w < \tilde{A}(0)).
\end{aligned}
$$

In the Markovian setting, $P_N$ in (2.10) can be replaced by $P$, and $\tilde{\tau}_z^-$ and $\tilde{M}(\tilde{\tau}_{-z}^-)$ do not depend on $\tilde{A}(0)$ since $\tilde{Y}(0+) = 0$. So we arrive at (2.8). Note that (2.10) is less informative than (2.9), equivalently, (2.3). Namely, we need (2.9) with $t = \infty$ to get (2.9) from (2.10).

We next sum (2.3) with $\varphi(y) = \mathbb{1}(y \leq t)$ over all $i, j, k$, which yields

COROLLARY 2.1. *If $v(i) = -1$ for $i \in S^-$, we have, for $t, x, z > 0$,*

$$P(\tau_0^+ \leq t, -Y(\tau_0^+-) \leq z, 0 < Y(\tau_0^+) \leq x)$$

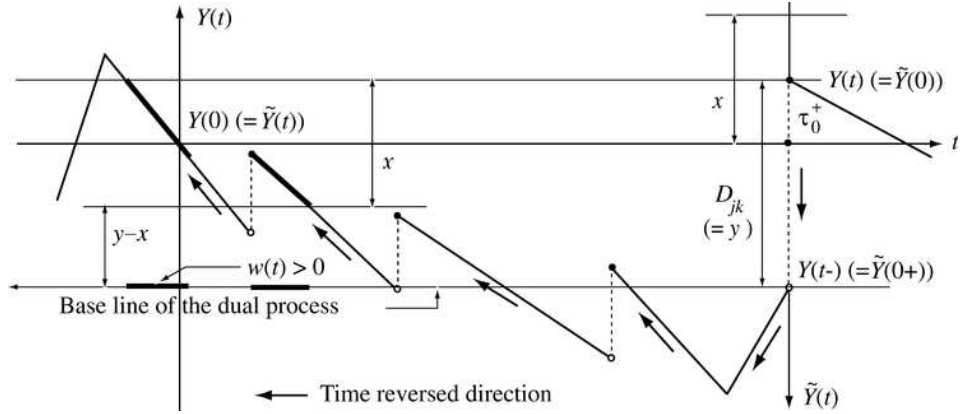

FIG. 1. *The additive process $Y(t)$ and its dual $\tilde{Y}(t)$.*

(2.11)
$$= \lambda \int_0^z P_N(\tilde{\tau}_{-w}^- \leq t, 0 < \tilde{A}(0) - w \leq x)\,dw,$$

where $v(i)$ can be an arbitrary positive number for $i \in S^+$.

REMARK 2.2. If $S^+ = \varnothing$ and $E(Y(1)) < 0$, then (2.11) with $t = \infty$ agrees with the well-known formula for the severity in a risk process [see, e.g., Asmussen (2000)], since $P(\tilde{\tau}_{-w}^- < \infty) = 1$ and $P_N(\tilde{A}(0) \leq x) = P_N(A(0) \leq x)$. On the other hand, if $x = z = \infty$, it generalizes Proposition 2.3 in Chapter IV of Asmussen (2000), which is obtained for the case that $N$ is Poisson.

Lemma 2.1 also generalizes Proposition 3.1 and Theorem 3.1 of Asmussen (1991), which assumes that $v(i) \equiv -1$. We here note another aspect of Lemma 2.1.

COROLLARY 2.2. *If there are no jumps, that is, $F_{ij}(0) = 1$ for all $i, j \in S$, then we have, for $i \in S^-$ and $k \in S^+$,*

(2.12) $\quad -v(i)P(M(0) = i, M(\tau_0^+) = k) = v(k)P(\tilde{M}(0) = k, \tilde{M}(\tilde{\tau}_0^-) = i).$

PROOF. Since $Y(\tau_0^+) = 0$ with probability one, (2.4) with $\varphi \equiv 1$ yields
$$-v(i)P(M(0) = i, M(\sigma^-(\tau_0^+)-) = j, M(\tau_0^+) = k)$$
$$= \lambda \int_0^\infty dw\, P_N(\tilde{M}(0) = k, \tilde{M}(0+) = j, \tilde{M}(\tilde{\tau}_{-w}^-) = i, w < \tilde{T}).$$

Summing over all possible $j$ and using the fact that $\tilde{T} = -v(k)\tilde{\sigma}^-(0)$ and $\tilde{M}(\tilde{\tau}_{-w}^-) \circ \theta_{-w/v(k)} = \tilde{M}(\tilde{\tau}_0^-)$ on $\tilde{M}(0) = k$, we have
$$-v(i)P(M(0) = i, M(\tau_0^+) = k)$$
$$= \lambda \int_0^\infty dw\, P_N(\tilde{M}(0) = k, \tilde{M}(\tilde{\tau}_{-w}^-) = i, w < \tilde{T})$$
$$= \lambda E_N\left(\int_0^{\tilde{T}} \mathbb{1}[\tilde{M}(-w/v(k)) = k, \tilde{M}(\tilde{\tau}_0^-) \circ \theta_{w/v(k)} = i]\,dw\right).$$

Since $\tilde{T} = v(k)\sigma^+(0)$, we get (2.12) by the Palm inversion formula. □

Equation (2.12) explains how upward hitting is interpreted in the dual additive process when there are no jumps. Suppose $-\infty < E(Y(1)) < 0$. Then (2.12) also yields
$$P(\tau_0^+ < \infty) = \frac{a^+}{a^-},$$
where $a^+ = \sum_{i \in S^+} v(i)P(M(0) = i)$ and $a^- = -\sum_{i \in S^-} v(i)P(M(0) = i)$.

**3. The Markov additive process.** From now on we specialize $M(t)$ to be a continuous-time Markov chain with finite state space $S$ and assume that $A(t)$ only depend on $M(t-)$ and $M(t)$. In this setting $(M(t), Y(t))$ is called a Markov additive process. In this section we derive the ascending ladder height distribution in a closed form for this additive process.

Let $C_{ij}$ for $i \ne j$ and $D_{ij}$ for all $i, j \in S$ be nonnegative numbers and let

$$C_{ii} = -\left(\sum_{j \ne i} C_{ij} + \sum_{j \in S} D_{ij}\right), \qquad i \in S.$$

For convenience, let $c(i) = -C_{ii}$. Let $C = \{C_{ij}\}$ and $D = \{D_{ij}\}$ be $S \times S$ matrices. We take $C + D$ for the rate matrix of Markov chain $M(t)$ and assume that $C + D$ is irreducible. Since $S$ is finite, this implies that the Markov chain has a unique stationary distribution, which is represented as a row vector $\boldsymbol{\pi} \equiv \{\pi_i\}$.

Note that the Markov chain $M(t)$ may include transitions that do not change their states, since $D_{ii}$ may be positive. Of course, these transitions are irrelevant for the sample path of $M(t)$, but they are convenient to include Markov modulated Poisson arrivals for the jumps in our formulation. So we let $N$ be the point process generated by all the transition instants. Let $F_{ij}^{\mathrm{D}}$ be the jump size distribution associated with transitions due to $D_{ij}$, where $F_{ij}^{\mathrm{D}}(0) = 0$, while no jumps are associated with transitions due to $C_{ij}$. Then, the conditional distribution $F_{ij}$ of $A(t)$, given that $M(t-) = i$ and $M(t) = j$, is obtained as

$$F_{ij}(x) = \begin{cases} F_{ii}^{\mathrm{D}}(x), & i = j, \\ \dfrac{C_{ij}}{C_{ij} + D_{ij}} \delta(x) + \dfrac{D_{ij}}{C_{ij} + D_{ij}} F_{ij}^{\mathrm{D}}(x), & i \ne j, \end{cases}$$

for $x \ge 0$, where $\delta(x)$ is the Dirac distribution which has a unit mass at the origin. In this way the additive process $Y(t)$ is completely specified for the rate function $v$.

Define $C_{ij}(x)$ and $D_{ij}(x)$ as

$$C_{ij}(x) = C_{ij}\delta(x), \qquad D_{ij}(x) = D_{ij}F_{ij}^{\mathrm{D}}(x), \qquad i, j \in S.$$

We denote the matrices which have these entries by $C(x)$ and $D(x)$, respectively. For convenience, we also define matrix $U(x)$ whose $ij$-entry is defined as

(3.1) $$U_{ij}(x) = \mathbb{1}(i \ne j)C_{ij}(x) + D_{ij}(x), \qquad i, j \in S, x \ge 0.$$

Clearly, the additive process $Y(t)$ is specified by $C$, $D(x)$ and the rate function $v$.

We now evaluate (2.3) and (2.4) in the Markov setting. To this end, we first consider the dual Markov chain $\tilde{M}(t)$ for $M(t)$ and the dual additive process $\tilde{Y}(t)$, defined in Section 2. Define $S \times S$ matrices:

$$\tilde{C}(x) = \Delta_{\boldsymbol{\pi}}^{-1} C(x)' \Delta_{\boldsymbol{\pi}}, \qquad \tilde{D}(x) = \Delta_{\boldsymbol{\pi}}^{-1} D(x)' \Delta_{\boldsymbol{\pi}},$$
$$\tilde{U}(x) = \Delta_{\boldsymbol{\pi}}^{-1} U(x)' \Delta_{\boldsymbol{\pi}},$$

where $\Delta_{\boldsymbol{\pi}}$ is the diagonal matrix whose $i$th entry is $\pi(i)$ and "$'$" stands for transpose. Obviously, we see that $\tilde{M}(t)$ has transition rate matrix $\tilde{C} + \tilde{D}$, where $\tilde{C} = \tilde{C}(\infty)$ and $\tilde{D} = \tilde{D}(\infty)$. $\tilde{Y}(t)$ and $Y(t)$ have the same rate function $v$ for the continuous movements. Thus, $\tilde{Y}(t)$ is the Markov additive process, specified by $\tilde{C}(x)$, $\tilde{D}(x)$ and $v$.

We next modify Lemma 2.1 in a convenient form for the Markov additive process. Let $\tilde{F}_{jk} = F_{kj}$. Then, from (2.3) and (2.4), we have, for $i \in S^-$, $j, k \in S$ and $x > 0$,

(3.2)
$$-v(i) E(M(0) = i, \mathbf{M}^-(\tau_0^+) = (j,k), 0 < Y(\tau_0^+) \leq x)$$
$$= \lambda \int_0^x dw \int_w^\infty \tilde{F}_{kj}(dy) P_N(\tilde{\mathbf{M}}^+(0) = (k,j), \tilde{M}(\tilde{\tau}_{-(y-w)}^-) = i),$$

and, for $i \in S^-$, $k \in S^+$,

(3.3)
$$-v(i) E(M(0) = i, \mathbf{M}^-(\tau_0^+) = (j,k), Y(\tau_0^+) = 0)$$
$$= \lambda \int_0^\infty dw \int_0^\infty \tilde{F}_{kj}(dy) P_N(\tilde{\mathbf{M}}^+(0) = (k,j), \tilde{M}(\tilde{\tau}_{-w}^-) = i,$$
$$y < w < \tilde{T} + y).$$

LEMMA 3.1. *For* $i \in S^-$, $j, k \in S$ *and* $x \geq 0$,

(3.4)
$$P(\mathbf{M}^-(\tau_0^+) = (j,k), 0 < Y(\tau_0^+) \leq x | M(0) = i)$$
$$= \frac{\pi_k}{\pi_i} \int_0^x dw \int_w^\infty \tilde{D}_{kj}(dy) P_N(\tilde{M}(\tilde{\tau}_{-(y-w)}^-) = i | \tilde{M}(0+) = j) \frac{-1}{v(i)}$$

*and, for* $i \in S^-$, $j \in S$ *and* $k \in S^+$,

(3.5)
$$P(\mathbf{M}^-(\tau_0^+) = (j,k), Y(\tau_0^+) = 0 | M(0) = i)$$
$$= \frac{\pi_k}{\pi_i} \int_0^\infty e^{-c(k)w/v(k)} dw \int_0^\infty \tilde{U}_{kj}(dy)$$
$$\times P_N(\tilde{M}(\tilde{\tau}_{-(w+y)}^-) = i | \tilde{M}(0+) = j) \frac{-1}{v(i)}.$$

PROOF. We only need to evaluate the Palm probabilities in (3.2) and (3.3). Since $\tilde{M}(t)$ is the Markov chain with rate matrix $\tilde{C} + \tilde{D}$, it is easy to see

that

$$\begin{aligned}
\lambda P_N(\tilde{M}(0) &= k, \tilde{M}(0+) = j, \tilde{M}(\tilde{\tau}^-_{-(y-w)}) = i) \\
&= P_N(\tilde{M}(\tilde{\tau}^-_{-(y-w)}) = i | \tilde{M}(0+) = j)(\mathbb{1}(k \neq j)\tilde{C}_{kj} + \tilde{D}_{kj})\pi_k.
\end{aligned} \quad (3.6)$$

Hence, (3.2) and (3.1) yields (3.4). Since $T$ of Lemma 2.1 is exponentially distributed with mean $-c(k)/v(k)$, the right-hand side of (3.3) becomes

$$\lambda \int_0^\infty \tilde{F}_{kj}(dy) \int_0^\infty dw \, P_N(\tilde{\mathbf{M}}^+(0) = (k,j), \tilde{M}(\tilde{\tau}^-_{-w}) = i, y < w < y + \tilde{T})$$

$$= \pi_k \int_0^\infty \tilde{U}_{kj}(dy) \int_0^\infty e^{-c(k)u/v(k)} \frac{-c(k)}{v(k)} \, du$$

$$\times \int_y^{y+u} dw \, P_N(\tilde{M}(\tilde{\tau}^-_{-w}) = i | \tilde{M}(0+) = j)$$

$$= \pi_k \int_0^\infty \tilde{U}_{kj}(dy) \int_y^\infty dw \int_{w-y}^\infty e^{-c(k)u/v(k)} \frac{-c(k)}{v(k)} \, du$$

$$\times P_N(\tilde{M}(\tilde{\tau}^-_{-w}) = i | \tilde{M}(0+) = j)$$

$$= \pi_k \int_0^\infty \tilde{U}_{kj}(dy) \int_y^\infty e^{-c(k)(w-y)/v(k)} \, dw$$

$$\times P_N(\tilde{M}(\tilde{\tau}^-_{-w}) = i | \tilde{M}(0+) = j).$$

Hence, (3.3) yields (3.5). □

In (3.4) and (3.5), it remains to evaluate

$$(3.7) \quad \tilde{R}_{ij}(x) \equiv P(\tilde{M}(\tilde{\tau}^-_{-x}) = j | \tilde{M}(0) = i), \qquad x \geq 0, \, i \in S, j \in S^-.$$

By $\tilde{R}^{--}(x)$ and $\tilde{R}^{+-}(x)$, we denote the $S^- \times S^-$ and $S^+ \times S^-$ matrices whose $ij$ entry is given by (3.7), respectively. We shall use this convention widely. Namely, for $S \times S$ matrix $A$, its submatrices of sizes $S^- \times S^-$, $S^- \times S^+$, $S^+ \times S^-$ and $S^+ \times S^+$ are denoted by $A^{--}$, $A^{-+}$, $A^{+-}$ and $A^{++}$, respectively. Similarly, for an $S$-vector $\mathbf{x}$, $\mathbf{x}^-$ and $\mathbf{x}^+$ are the subvectors whose entries are indexed by $S^-$ and $S^+$, respectively. Note that $\tilde{R}^{--}(x)$ is substochastic, that is, a nonnegative matrix satisfying $\tilde{R}^{--}(x)\mathbf{e}^- \leq \mathbf{e}^-$, where $\mathbf{e}$ is the column vector whose entries are all 1.

We now refer to recent results on the matrix exponential form in the Markov additive process in Miyazawa and Takada (2002) and Takada (2001). The results there are obtained for the forward process $\{(Y(t), M(t))\}$ that is right-continuous. On the other hand, we here need the corresponding results for the dual process $\{(\tilde{Y}(t), \tilde{M}(t))\}$ that is left-continuous. However, this left-continuity is irrelevant in the definition (3.7), since $\tilde{Y}(t)$ is skip free in the downward direction.

This skip free property also implies that the ladder height indexed process $\{\tilde{M}(\tilde{\tau}^-_{-x}); x \geq 0\}$ is Markov. Hence, $\tilde{R}^{--}(x) = e^{x\tilde{Q}^{(-)}}$ for some subrate matrix $\tilde{Q}^{(-)}$, that is, $\tilde{Q}^{(-)}\mathbf{e}^- \leq \mathbf{0}^-$. It is also easy to see, at least intuitively, that

$$-\frac{d}{dx}\begin{pmatrix}\tilde{R}^{--}(x)\\ \tilde{R}^{+-}(x)\end{pmatrix} = \Delta_\mathbf{v}^{-1}\int_0^\infty (\tilde{C}(du) + \tilde{D}(du))\begin{pmatrix}\tilde{R}^{--}(u+x)\\ \tilde{R}^{+-}(u+x)\end{pmatrix}.$$

Then, the following results would be intuitively reasonable. Their formal proofs not for $\tilde{Y}(t)$, but for the forward process $Y(t)$ are given in Theorem 3.2 of Miyazawa and Takada (2002) and Theorem 4.1 of Takada (2001). Clearly, the results are immediately transferred from $Y(t)$ to $\tilde{Y}(t)$. In what follows, a square matrix is said to be ML if its off-diagonal entries are nonnegative.

LEMMA 3.2. *There exist a subrate matrix $\tilde{Q}^{(-)}$ and a substochastic matrix $\tilde{R}^{+-}$, such that*

$$(3.8) \qquad \begin{pmatrix}\tilde{R}^{--}(x)\\ \tilde{R}^{+-}(x)\end{pmatrix} = \begin{pmatrix}I^{--}\\ \tilde{R}^{+-}\end{pmatrix}e^{x\tilde{Q}^{(-)}}, \qquad x \geq 0,$$

*where $\tilde{Q}^{(-)}$ and $\tilde{R}^{+-}$ are the minimal ML and nonnegative solutions, respectively, of the following equation:*

$$(3.9) \qquad -\begin{pmatrix}I^{--}\\ \tilde{R}^{+-}\end{pmatrix}\tilde{Q}^{(-)} = \Delta_\mathbf{v}^{-1}\int_0^\infty (\tilde{C}(du) + \tilde{D}(du))\begin{pmatrix}I^{--}\\ \tilde{R}^{+-}\end{pmatrix}e^{u\tilde{Q}^{(-)}}.$$

*Furthermore, $\tilde{Q}^{(-)}$ is a rate matrix, if and only if*

$$(3.10) \qquad E(Y(1))(=E(\tilde{Y}(1))) = \boldsymbol{\pi}\Delta_\mathbf{v}\mathbf{e} + \boldsymbol{\pi}\int_0^\infty xD(dx)\mathbf{e} \leq 0.$$

We next represent $\tilde{Q}^{(-)}$ and $\tilde{R}^{+-}$ using the original additive process. Let $K^{(-)}$ and $L^{-+}$ be the minimal solutions of the following equation, such that $K^{(-)}$ and $L^{-+}$ are ML and nonnegative matrices, respectively:

$$(3.11) \quad -K^{(-)}(I^{--}, L^{-+}) = \int_0^\infty e^{uK^{(-)}}(I^{--}, L^{-+})(C(du) + D(du))\Delta_\mathbf{v}^{-1}.$$

Then $\tilde{Q}^{(-)}$ and $\tilde{R}^{+-}$ are obtained as

$$(3.12) \quad \tilde{Q}^{(-)} = (\Delta_{\boldsymbol{\pi}}^{--})^{-1}(K^{(-)})'\Delta_{\boldsymbol{\pi}}^{--}, \qquad \tilde{R}^{+-} = (\Delta_{\boldsymbol{\pi}}^{++})^{-1}(L^{-+})'\Delta_{\boldsymbol{\pi}}^{--}.$$

COROLLARY 3.1. *If (3.10) is satisfied, then*

$$(3.13) \qquad \boldsymbol{\pi}^- K^{(-)} = \mathbf{0}, \qquad \boldsymbol{\pi}^- L^{-+} = \boldsymbol{\pi}^+,$$

*and $K^{(-)}$ has the right eigenvector $\mathbf{k}^-$ corresponding to $\boldsymbol{\pi}^-$. Normalize $\mathbf{k}^-$ in such a way that $\boldsymbol{\pi}^-\mathbf{k}^- = 1$. Then, $\mathbf{k}^-$ is unique and positive and $\mathbf{k}^-\boldsymbol{\pi}^- - K^{(-)}$ is nonsingular.*

PROOF. Since $\tilde{Q}^{(-)}\mathbf{e}^- = \mathbf{0}^-$ and $\tilde{R}^{+-}\mathbf{e}^- = \mathbf{e}^+$ by (3.10), (3.12) implies (3.13). $(\mathbf{k}^-)'\Delta_{\boldsymbol{\pi}}^{--}$ is the stationary distribution of $\tilde{Q}^{(-)}$, since $(\mathbf{k}^-)'\Delta_{\boldsymbol{\pi}}^{--}\mathbf{e} = \boldsymbol{\pi}^-\mathbf{k}^- = 1$. Hence, $\mathbf{k}^-$ is unique and positive. Assume that

$$(3.14) \qquad (\mathbf{k}^-\boldsymbol{\pi}^- - K^{(-)})\mathbf{x} = \mathbf{0},$$

for an $S^-$-dimensional vector $\mathbf{x}$. Pre-multiplying $K^{(-)}$ to (3.14), we get $K^{(-)}(K^{(-)}\mathbf{x}) = \mathbf{0}$. Hence, $K^{(-)}\mathbf{x}$ must be $a\mathbf{k}^-$ for some constant $a$. Since $\boldsymbol{\pi}^- K^{(-)} = \mathbf{0}$ and $\boldsymbol{\pi}^-\mathbf{k}^- = 1$, we arrive at $a = 0$, which concludes $\mathbf{x} = \mathbf{0}$. So the matrix $\mathbf{k}^-\boldsymbol{\pi}^- - K^{(-)}$ is nonsingular. $\square$

We remark that similar eigenvectors have been reported for the case that $S^+ = \varnothing$ in Section 2 of Chapter VI in Asmussen (2000). Let us compute the ascending ladder height distribution. For this computation we need one more lemma.

LEMMA 3.3. *Let $\Delta_{\mathbf{c}}$ be the diagonal matrix whose $j$th diagonal entry is $c(j) = -C_{jj}$, then, for $j \in S^+$ and $i \in S^-$,*

$$(3.15) \quad \int_0^\infty e^{wc(j)/v(j)}\,dw\left[\int_0^\infty \tilde{U}(dy)\begin{pmatrix}I^{--}\\ \tilde{R}^{+-}\end{pmatrix}e^{(y+w)\tilde{Q}^{(-)}}\right]_{ji} = [\Delta_{\mathbf{v}}^{++}\tilde{R}^{+-}]_{ji}.$$

PROOF. By (3.9), the left-hand side of (3.15) is

$$\int_0^\infty e^{wc(j)/v(j)}\,dw\left[\int_0^\infty (\tilde{C}(dy) + \tilde{D}(dy) - \delta(dy)(-\Delta_{\mathbf{c}}))\begin{pmatrix}I^{--}\\ \tilde{R}^{+-}\end{pmatrix}e^{(y+w)\tilde{Q}^{(-)}}\right]_{ji}$$

$$= -\int_0^\infty e^{wc(j)/v(j)}\,dw\left[\left(\Delta_{\mathbf{v}}\begin{pmatrix}I^{--}\\ \tilde{R}^{+-}\end{pmatrix}\tilde{Q}^{(-)} + (-\Delta_{\mathbf{c}})\begin{pmatrix}I^{--}\\ \tilde{R}^{+-}\end{pmatrix}\right)e^{w\tilde{Q}^{(-)}}\right]_{ji}$$

$$= [\Delta_{\mathbf{v}}^{++}\tilde{R}^{+-}]_{ji}. \qquad \square$$

THEOREM 3.1. *For $i \in S^-$, $j \in S$ and $x \geq 0$,*

$$P(M(\tau_0^+) = j, Y(\tau_0^+) \leq x | M(0) = i)$$

$$(3.16) \qquad = \frac{1}{-v(i)}\left(\int_0^x dw\left[\int_w^\infty e^{(y-w)K^{(-)}}(I^{--}, L^{-+})D(dy)\right]_{ij}\right.$$

$$\left. + \mathbb{1}(j \in S^+)L_{ij}^{-+}v(j)\right).$$

PROOF. Let $J_{ij}(x) = P(M(\tau_0^+) = j, Y(\tau_0^+) \leq x | M(0) = i)$. We first consider the case that $Y(\tau_0^+) = 0$. From (3.5) and Lemma 3.3, we have

$$J_{ij}(0) = \frac{\pi_j}{\pi_i}v(j)\tilde{R}_{ji}^{+-}\frac{-1}{v(i)}.$$

Using $\pi_j \tilde{R}_{ji}^{+-} = \pi_i L_{ij}^{-+}$, this yields (3.16) for $x = 0$. We next consider the case that $Y(\tau_0^+) > 0$. From Lemmas 3.1 and 3.2, we have

$$J_{ij}(x) - J_{ij}(0) = \frac{\pi_j}{\pi_i} \int_0^x dw \left[ \int_w^\infty \tilde{D}(dy) \begin{pmatrix} I^{--} \\ \tilde{R}^{+-} \end{pmatrix} e^{(y-w)\tilde{Q}^{(-)}} \right]_{ji} \frac{-1}{v(i)}.$$

Thus, we get (3.16) by converting to the notation of the forward processes. □

**4. The hitting probability.** We now consider the hitting probability at an upper level $x > 0$. Since $(X(t), Y(t))$ is a Markov additive process with a real additive component, we can obtain the hitting probabilities as a Markov renewal function using the appropriate ladder height distribution as a semi-Markov kernel. However, we can not use (3.16) of Theorem 3.1 as this kernel because the additive component may increase continuously. To have an appropriate kernel, we consider the ladder height at $\sigma^+(\tau_0^+)$, that is, the first transition instant of $M(t)$ after crossing level 0. For convenience, define a random variable $\xi_x^+$ for $x \geq 0$ as

$$\xi_x^+ = \sigma^+(\tau_x^+).$$

Clearly, $\xi_x^+$ is a stopping time with respect to the additive process $\{(Y(t), M(t))\}$. Define $S \times S$ matrices $H(x)$ and $\overline{G}(x)$ by

$$H_{ij}(x) = P(M(\xi_0^+) = j, Y(\xi_0^+) \leq x | M(0) = i),$$
$$\overline{G}_{ij}(x) = P(M(\tau_0^+) = j, Y(\xi_0^+) \geq x | M(0) = i), \qquad i, j \in S, x \geq 0.$$

Define the hitting probability matrix $\Psi(x)$ by

$$\Psi_{ij}(x) = P(M(\tau_x^+) = j | M(0) = i).$$

Then, $\Psi(x)$ is obtained as the unique solution of the Markov renewal equation,

(4.1) $$\Psi(x) = \overline{G}(x) + H * \Psi(x), \qquad x \geq 0,$$

where the convolution $A * B(x)$ is defined for a matrix nondecreasing function $A(x)$ and a matrix function $B(x)$ as

$$[A * B(x)]_{ij} = \sum_k \int_0^\infty A_{ik}(dy) B_{kj}(x - y).$$

The Markov renewal equation (4.1) will be the key for our arguments to obtain an asymptotic behavior of the hitting probability. Thus, what we all need is to get $H(x)$ and $\overline{G}(x)$. We first consider the case that $Y(\tau_0^+) = 0$.

LEMMA 4.1. *For $i \in S^-$, $j \in S$ and $x \geq 0$,*

$$(4.2) \quad \begin{aligned} &P(M(\sigma^+(\tau_0^+)) = j, Y(\tau_0^+) = 0, Y(\sigma^+(\tau_0^+)) \leq x | M(0) = i) \\ &= \frac{\pi_j}{\pi_i} \sum_{k \in S^+} \int_0^x \tilde{U}_{jk}(dy)(1 - e^{-c(k)(x-y)/v(k)}) \frac{v(k)}{c(k)} \tilde{R}_{ki}^{+-} \frac{-1}{v(i)}. \end{aligned}$$

PROOF. We need the following version of (3.3), for $i \in S^-$, $j, l \in S$ and $k \in S^+$:

$$(4.3) \quad \begin{aligned} -v(i)P(M(0) = i, \mathbf{M}^-(\tau_0^+) = (l, k), M(\sigma^+(\tau_0^+)) = j, \\ Y(\tau_0^+) = 0, Y(\sigma^+(\tau_0^+)) \leq x) \\ = \lambda \int_0^\infty dw \int_0^\infty \tilde{F}_{jk}(dy) P_N(\tilde{M}(\tilde{\sigma}^-(0)) = j, \\ \tilde{\mathbf{M}}^+(0) = (k, l), \tilde{M}(\tilde{\tau}_{-w}^-) = i, \\ y < w < \tilde{T} + y, \tilde{T} + \tilde{T}' + y \leq x + w), \end{aligned}$$

where $\tilde{T}' = -\tilde{Y}(\sigma^-(0)) + \tilde{Y}(\sigma^-(0)+)$. Since this can be proved in exactly the same way as (3.3), we omit its proof. Thus, the right-hand side of (4.3) divided by $-v(i)$ becomes

$$-\frac{\lambda}{v(i)} \int_0^\infty dw \int_0^\infty \tilde{F}_{jk}(dz) \int_0^\infty \frac{c(k)}{v(k)} e^{-c(k)u/v(k)} du \int_0^\infty \tilde{F}_{kl}(dy)$$
$$\times P_N(\tilde{M}(\tilde{\sigma}^-(0)) = j, \tilde{\mathbf{M}}^+(0) = (k, l), \tilde{M}(\tilde{\tau}_{-w}^-) = i,$$
$$y \leq w < u + y, u + y + z \leq x + w)$$

(changing variable $w$ to $w - y$)

$$= -\frac{\lambda}{v(i)} \int_0^\infty dw \int_0^\infty \tilde{F}_{kl}(dy) \int_0^\infty \frac{c(k)}{v(k)} e^{-c(k)u/v(k)} du \int_0^\infty \tilde{F}_{jk}(dz)$$
$$\times P_N(\tilde{M}(\tilde{\sigma}^-(0)) = j, \tilde{\mathbf{M}}^+(0) = (k, l), \tilde{M}(\tilde{\tau}_{-(y+w)}^-) = i,$$
$$0 \leq w < u, u + z \leq x + w)$$
$$= -\frac{\lambda}{v(i)} \int_0^\infty dw \int_0^\infty \tilde{F}_{kl}(dy) \int_w^{x+w} \frac{c(k)}{v(k)} e^{-c(k)u/v(k)} du \int_0^{x+w-u} \tilde{F}_{jk}(dz)$$
$$\times P_N(\tilde{M}(\tilde{\sigma}^-(0)) = j, \tilde{\mathbf{M}}^+(0) = (k, l), \tilde{M}(\tilde{\tau}_{-(y+w)}^-) = i),$$

where we have used the following fact. Conditionally, given that $(\tilde{\sigma}^-(0)) = j$, $\tilde{M}(0) = k, \tilde{M}(0+) = l, \tilde{M}(\tilde{\tau}_{-w}^-) = i,$

$$-\tilde{Y}(0) \simeq \tilde{F}_{kl},$$
$$\tilde{T} = -\tilde{Y}(\tilde{\sigma}^-(0)+) + \tilde{Y}(0) \simeq \text{Exp}(c(k)/v(k)),$$
$$\tilde{T}' = -\tilde{Y}(\tilde{\sigma}^-(0)) + \tilde{Y}(\tilde{\sigma}^-(0)+) \simeq \tilde{F}_{jk},$$

where $\simeq$ stands for equality in distribution and $\text{Exp}(a)$ denotes the exponential distribution with mean $1/a$. Hence, using a decomposition similarly to (3.6) and Lemma 3.2, the left-hand side of (4.2) is computed as

$$\frac{\pi_j}{\pi_i} \int_0^\infty \left( \int_0^\infty dw \int_w^{x+w} du \int_0^{x+w-u} (\mathbb{1}(j \neq k)\tilde{C}_{jk} + \tilde{D}_{jk})\tilde{F}_{jk}(dz) \right.$$

$$\times \frac{c(k)}{v(k)} e^{-c(k)u/v(k)} \Big)$$

$$\times \frac{\mathbb{1}(k \neq l)\tilde{C}_{kl} + \tilde{D}_{kl}}{c(k)} \tilde{F}_{kl}(dy) \left[ \begin{pmatrix} I^{--} \\ \tilde{R}^{+-} \end{pmatrix} e^{(y+w)\tilde{Q}^{(-)}} \right]_{li} \frac{-1}{v(i)}$$

$$= \frac{\pi_j}{\pi_i} \int_0^\infty \left( \int_0^\infty e^{-c(k)w/v(k)} \, dw \int_0^x du \right.$$

$$\times \int_0^{x-u} \tilde{U}_{jk}(dz) \frac{c(k)}{v(k)} e^{-c(k)u/v(k)} \Big) \frac{1}{c(k)} \tilde{U}_{kl}(dy)$$

$$\times \left[ \begin{pmatrix} I^{--} \\ \tilde{R}^{+-} \end{pmatrix} e^{(y+w)\tilde{Q}^{(-)}} \right]_{li} \frac{-1}{v(i)}$$

$$= \frac{\pi_j}{\pi_i} \int_0^\infty \left( \int_0^\infty e^{-c(k)w/v(k)} \, dw \int_0^x \tilde{U}_{jk}(dz) \right.$$

$$\times (1 - e^{-c(k)(x-z)/v(k)}) \Big) \frac{-1}{-c(k)} \tilde{U}_{kl}(dy)$$

$$\times \left[ \begin{pmatrix} I^{--} \\ \tilde{R}^{+-} \end{pmatrix} e^{(y+w)\tilde{Q}^{(-)}} \right]_{li} \frac{-1}{v(i)}.$$

Summing this over all $l \in S$ and $k \in S^+$ and applying Lemma 3.3 yield (4.2). $\square$

THEOREM 4.1. *For $i, j \in S$ and $x \geq 0$, $H_{ij}(x)$ is given by*

(4.4)
$$H_{ij}(x) = \frac{1}{-v(i)} \left( \int_0^x dw \left[ \int_w^\infty e^{(y-w)K^{(-)}} (I^{--}, L^{-+}) D(dy) \right]_{ij} \right.$$

$$+ \sum_{k \in S^+} \int_0^x L_{ik}^{-+} (1 - e^{-c(k)(x-y)/v(k)}) \frac{v(k)}{c(k)} U_{kj}(dy) \right),$$

$$i \in S^-$$

$$H_{ij}(x) = \int_0^x (1 - e^{-c(i)(x-y)/v(i)}) \frac{1}{c(i)} U_{ij}(dy), \qquad i \in S^+.$$

PROOF. We first consider the case that $i \in S^-$. From Lemma 4.1, we get $H_{ij}(0)$ of (4.4). On the other hand, from Theorem 3.1,

$$P(M(\tau_0^+) = j, 0 < Y(\tau_0^+) \le x | M(0) = i)$$
$$= \frac{1}{-v(i)} \int_0^x dw \left[ \int_w^\infty e^{(y-w)K^{(-)}} (I^{--}, L^{-+}) D(dy) \right]_{ij}.$$

Thus, we get (4.4) converting to the notation of the forward processes. For $i \in S^+$, (4.4) is immediate, since $\xi_0^+$ is the first transition instant after time 0 in this case. □

For $\overline{G}$, we can get the following expression similarly to Theorem 4.1, using the fact that the first hitting over level $x$ is attained either continuously or by a jump.

COROLLARY 4.1. *For $i \in S^+, j \in S$,*

(4.5)
$$\overline{G}_{ij}(x) = \mathbb{1}(i=j)e^{-xc(i)/v(i)}$$
$$+ \int_0^x \frac{c(i)}{v(i)} e^{-c(i)y/v(i)} \, dy \int_{x-y}^\infty \frac{1}{c(i)} U_{ij}(dw), \qquad x \ge 0,$$

*and, for $i \in S^-, j \in S$,*

(4.6)
$$\overline{G}_{ij}(x) = \frac{-1}{v(i)} \left[ \int_x^\infty dw \left[ \int_w^\infty e^{(y-w)K^{(-)}} (I^{--}, L^{-+}) D(dy) \right]_{ij} \right.$$
$$+ L_{ij}^{-+} v(j) e^{xC_{jj}/v(j)}$$
$$\left. + \sum_{k \in S^+} L_{ik}^{-+} \frac{v(k)}{c(k)} \int_0^x \int_{x-y}^\infty U_{kj}(dw) \frac{c(k)}{v(k)} e^{-yc(k)/v(k)} \, dy \right],$$
$$x \ge 0.$$

**5. Asymptotic behavior of the hitting probability.** In this section, we study the asymptotic behavior of the hitting probabilities as the hitting level goes to infinity. Throughout this section we assume a negative drift, that is,

(5.1) $$E(Y(1)) = \pi \Delta_{\mathbf{v}} \mathbf{e} + \pi \int_0^\infty D(du) \mathbf{e} < 0.$$

Under this condition, $Y(t) \to -\infty$ as time $t \to \infty$. Hence, the hitting probability $P(\tau_x^+ < \infty | M(0) = i)$ converges to zero as $x \to \infty$. Assume that the jump size distributions have light tails. Then, it is expected that the hitting probability decays exponentially fast, which is known as *the Cramér–Lundberg approximation* for the conventional risk model. Furthermore, a

Brownian component may be added [see Schmidli (1995)]. Instead of a Brownian component, we here have signed continuous movements.

The decay rates for the hitting probabilities have been extensively studied by the change of measure technique based on martingales [see, e.g., Asmussen (2000) and Rolski, Schmidli, Schmidt and Teugels (1999)], but this approach has not yet fully covered the Markov modulated model even for the case that $v(i) \equiv -1$ [see Section VI of Asmussen (2000)]. Here we use the Markov renewal theorem as in Miyazawa (2002), which considers the case that $v(i) \equiv -1$. This approach also uses a change of measure, but it is more straightforward.

Let us briefly introduce this Markov renewal approach, following Section 2 of Miyazawa (2002). Let $S \times S$ matrix $P(x) = \{P_{ij}(x)\}$ be a Markov renewal kernel that may be defective, that is, $P(x)\mathbf{e} \le \mathbf{e}$. Denote the moment generating function of $P(x)$ by $\hat{P}(\theta) \equiv \int_0^\infty e^{\theta x} P(dx)$. Since $\hat{P}(\theta)$ is a nonnegative and substochastic matrix, it has a positive eigenvalue $\gamma(\theta)$, such that the absolute values of all other eigenvalues are less than $\gamma(\theta)$ and the associated right and left eigenvectors are nonnegative [see, e.g., Seneta (1980)]. Denote these associated eigenvectors by $\boldsymbol{\nu}^{(\theta)}$ and $\mathbf{h}^{(\theta)}$, respectively. Suppose that an $S \times S$ matrix functions $\Phi(x)$ and $B(x)$ for $x \ge 0$ satisfies the Markov renewal equation,

$$(5.2) \qquad \Phi(x) = B(x) + P * \Phi(x), \qquad x \ge 0.$$

Then, Theorem 2.6 of Chapter X in Asmussen (1987) [see Miyazawa (2002) for the present context] reads as

LEMMA 5.1. *Suppose the following four conditions:* (5.a) $P(x)$ *has a single irreducible recurrent class that can be reached from any state in $S$ with probability one, and the return time to each state in the irreducible recurrent class has a nonarithmetic distribution,* (5.b) *there exists a positive $\alpha$ such that $\gamma(\alpha) = 1$,* (5.c) *each entry of $e^{\alpha x} B(x)$ is directly Riemann integrable,* (5.d) $\boldsymbol{\nu}^{(\alpha)} \hat{P}_{(1)}(\alpha) \mathbf{h}^{(\alpha)}$ *is finite. Then, we have*

$$(5.3) \qquad \lim_{x \to \infty} e^{\alpha x} \Phi(x) = \frac{1}{\boldsymbol{\nu}^{(\alpha)} \hat{P}_{(1)}(\alpha) \mathbf{h}^{(\alpha)}} \mathbf{h}^{(\alpha)} \boldsymbol{\nu}^{(\alpha)} \int_0^\infty e^{\alpha u} B(u) \, du,$$

*where* $\hat{P}_{(1)}(\alpha) = \frac{d}{d\theta} \hat{P}(\theta)|_{\theta = \alpha}$.

Let us put $P(x) = H(x)$, $\Phi(x) = \Psi(x)$ and $B(x) = \overline{G}(x)$. Then, (5.2) holds by (4.1). Clearly, condition (5.a) is satisfied by the irreducibility of $C + D$ and the exponential sojourn times of $M(t)$. We next compute the moment generating function $\hat{H}(\theta)$ of $H(x)$. From Theorem 4.1, the following results are obtained.

Lemma 5.2.

(5.4) $$\hat{H}(\theta) = I - \Delta_{\mathbf{v}}^{-1} T(\theta)(C + \hat{D}(\theta) + \theta \Delta_{\mathbf{v}}),$$

*where*

$$T(\theta) = \begin{bmatrix} (\theta I^{--} - K^{(-)})^{-1} & (\theta I^{--} - K^{(-)})^{-1} L^{-+} + L^{-+} \Delta_{\mathbf{v}}^{++} (\Delta_{\mathbf{c}-\theta\mathbf{v}}^{++})^{-1} \\ 0^{+-} & -\Delta_{\mathbf{v}}^{++} (\Delta_{\mathbf{c}-\theta\mathbf{v}}^{++})^{-1} \end{bmatrix}.$$

REMARK 5.1. Since $T(\theta)$ is invertible as we will see, (5.4) can be written as

(5.5) $$C + \hat{D}(\theta) + \theta \Delta_{\mathbf{v}} = T(\theta)^{-1} \Delta_{\mathbf{v}} (I - \hat{H}(\theta)).$$

This can be considered as a generator version of the Wiener–Hopf factorization [see, e.g., Arjas and Speed (1973)].

Since this lemma is just computations, we defer its proof to Appendix A. We notice that $\hat{H}(\theta)$ exists only for $\theta < \min\{c(i)/v(i); i \in S^+\}$. We shall use the following light tail assumption on $D(x)$ and related regularity assumption:

(i) $D(\theta)$ is finite for some $\theta > 0$,
(ii) $\sup\{D'_{ij}(\theta) < \infty; i,j \in S, \theta > 0\} = \infty$.

Condition (ii) can be relaxed, but it is sufficient for most applications.

Note that condition (5.1) implies that $\tilde{Q}^{(-)}$ is a rate matrix. Since $C + \hat{D}(\theta) + \theta \Delta_{\mathbf{v}}$ is an ML matrix for each $\theta > 0$, it has a real eigenvalue $\kappa(\theta)$ such that it dominates the real parts of all other eigenvalues, and the associated left and right eigenvectors are positive. Denote these eigenvectors by $\boldsymbol{\mu}^{(\theta)}$ and $\mathbf{h}^{(\theta)}$, respectively. Let $\mathbf{k}^-$ be the right eigenvector of $K^{(-)}$ for eigenvalue 0, where $\mathbf{k}^-$ is unique and positive (see Corollary 3.1). We normalize $\boldsymbol{\mu}^{(\theta)}$ and $\mathbf{h}^{(\theta)}$ so that

$$(\boldsymbol{\mu}^{(\theta)})^- \mathbf{k}^- = 1, \qquad \boldsymbol{\mu}^{(\theta)} \mathbf{h}^{(\theta)} = 1.$$

Since $\kappa(\theta) = \boldsymbol{\mu}^{(\theta)} (\theta \Delta_{\mathbf{v}} + C + \hat{D}(\theta)) \mathbf{h}^{(\theta)}$, we have

$$\kappa'(\theta) = \kappa(\theta)((\boldsymbol{\mu}^{(\theta)})' \mathbf{h}^{(\theta)} + \boldsymbol{\mu}^{(\theta)} (\mathbf{h}^{(\theta)})') + \boldsymbol{\mu}^{(\theta)} (\Delta_{\mathbf{v}} + \hat{D}'(\theta)) \mathbf{h}^{(\theta)}$$
$$= \boldsymbol{\mu}^{(\theta)} (\Delta_{\mathbf{v}} + \hat{D}'(\theta)) \mathbf{h}^{(\theta)},$$

where we have used the fact that $\boldsymbol{\mu}^{(\theta)} \mathbf{h}^{(\theta)} = 1$. This implies that $\kappa'(0) = E(Y(1)) < 0$ because $\boldsymbol{\mu}^{(0)} = \boldsymbol{\pi}$. Hence, from the fact that $\kappa(0) = 0$ and $\kappa(\theta)$ is a convex function [see Kingman (1961)], condition (ii) guarantees that $\kappa(\theta) = 0$ has a unique positive solution. Thus, condition (5.b) is satisfied. Denote this solution by $\alpha$.

To find $\boldsymbol{\nu}^{(\alpha)}$ for Lemma 5.1, we compute the inverse matrix of $T(\alpha)$:

$$T(\alpha)^{-1} = \begin{bmatrix} \alpha I^{--} - K^{(-)} & L^{-+}\Delta_{\mathbf{c}-\alpha\mathbf{v}}^{++}(\Delta_{\mathbf{v}}^{++})^{-1} + (\alpha I^{--} - K^{(-)})L^{-+} \\ 0^{+-} & -\Delta_{\mathbf{c}-\alpha\mathbf{v}}^{++}(\Delta_{\mathbf{v}}^{++})^{-1} \end{bmatrix}.$$

Define $\boldsymbol{\nu}^{(\alpha)}$ as

(5.6) $$\boldsymbol{\nu}^{(\alpha)} = -\boldsymbol{\mu}^{(\alpha)}T(\theta)^{-1}\Delta_{\mathbf{v}}.$$

Clearly, $\boldsymbol{\nu}^{(\alpha)}$ is the left invariant vector of $\hat{H}(\alpha)$. We show that $\boldsymbol{\nu}^{(\alpha)}$ is positive. To this end, we introduce a twisted Markov transition kernel for $H$. Define

$$H^\dagger(x) = \Delta_{\mathbf{h}^{(\alpha)}}^{-1} \int_0^x e^{\alpha u} H(du) \Delta_{\mathbf{h}^{(\alpha)}}, \qquad x \geq 0.$$

Then, by our choice of $\mathbf{h}^{(\alpha)}$, $H^\dagger(\infty)\mathbf{e} = \mathbf{e}$, that is, $H^\dagger(\infty)$ is a stochastic kernel. From (5.4), we have

$$\boldsymbol{\nu}^{(\alpha)}\Delta_{\mathbf{h}^{(\alpha)}}H^\dagger(\infty) = \boldsymbol{\nu}^{(\alpha)}\hat{H}(\alpha)\Delta_{\mathbf{h}^{(\alpha)}}$$
$$= \boldsymbol{\nu}^{(\alpha)}\Delta_{\mathbf{h}^{(\alpha)}} + \boldsymbol{\mu}^{(\alpha)}(\alpha\Delta_{\mathbf{v}} + C + \hat{D}(\alpha)) = \boldsymbol{\nu}^{(\alpha)}\Delta_{\mathbf{h}^{(\alpha)}}.$$

Hence, $\boldsymbol{\nu}^{(\alpha)}\Delta_{\mathbf{h}^{(\alpha)}}$ is the left invariant vector of stochastic kernel $H^\dagger(\infty)$. From (5.6) and the normalization of $\boldsymbol{\mu}^{(\alpha)}$, it follows that $(\boldsymbol{\nu}^{(\alpha)})^-(-\Delta_{\mathbf{v}}^{--})^{-1}\mathbf{k}^- = \alpha$, and $\mathbf{k}^-$ is positive by Corollary 3.1. So, $\boldsymbol{\nu}^{(\alpha)}$ must be nonnegative and $\boldsymbol{\nu}^{(\alpha)}\Delta_{\mathbf{h}^{(\alpha)}}$ is also nonnegative. Since $H^\dagger(\infty)$ is a finite stochastic matrix and irreducible, which follows from the irreducibility of $C + D$, $\boldsymbol{\nu}^{(\alpha)}\Delta_{\mathbf{h}^{(\alpha)}}$ must be a positive vector and unique up to a multiplicative constant. Hence, $\boldsymbol{\nu}^{(\alpha)}$ is a positive vector.

It is easy to see that $\mathbf{h}^{(\alpha)}$ is the right positive eigenvector of $\hat{H}(\alpha)$ for the eigenvalue 1. Thus, we are now in a position to apply Lemma 5.1, for $P = H$ and $B = \overline{G}$. We first compute the following integrals. For $i \in S^+$, (4.5) yields

(5.7)
$$\int_0^\infty e^{\alpha x}\overline{G}_{ij}(x)\,dx$$
$$= \delta_{ij}\int_0^\infty e^{\alpha x}e^{-c(i)/(v(i))x}\,dx$$
$$+ \frac{1}{v(i)}\int_0^\infty e^{\alpha x}\,dx \int_0^x e^{-c(i)/(v(i))y}\,dy \int_{x-y}^\infty U_{ij}(du)$$
$$= \frac{1}{\alpha}[\Delta_{\mathbf{c}-\alpha\mathbf{v}}^{-1}(C + \hat{D}(\alpha) + \alpha\Delta_{\mathbf{v}} - (C + D))]_{ij}.$$

Similarly, for $i \in S^-$, we have

$$\int_0^\infty e^{\alpha x}\overline{G}_{ij}(x)\,dx$$

$$= \frac{-1}{\alpha v(i)} \Big[ (\alpha I^{--} - K^{(-)})^{-1}(I^{--}, L^{-+})(C + \hat{D}(\alpha) + \alpha \Delta_{\mathbf{v}})$$

(5.8)
$$+ (0^{--}, L^{-+}) \Delta_{\mathbf{v}} \Delta_{\mathbf{c}-\alpha \mathbf{v}}^{-1}(C + \hat{D}(\alpha) + \alpha \Delta_{\mathbf{v}} - (C + D))$$

$$- \mathbf{k}^{-}\boldsymbol{\pi}^{-}(I^{--}, L^{-+}) \Big( \Delta_{\mathbf{v}} + \int_0^{\infty} y D(dy) \Big)$$

$$- (\mathbf{k}^{-}\boldsymbol{\pi}^{-} - K^{(-)})^{-1}(I^{--}, L^{-+})(C + D) \Big]_{ij}.$$

See Section A.2 for the detailed derivation of this formula. Thus, condition (5.c) is satisfied. Then, applying Lemma 5.1, we have the following asymptotics.

THEOREM 5.1. *Suppose the stability condition* (5.1) *and the light tail conditions* (i) *and* (ii)*. Then, the maximal eigenvalue* $\kappa(\theta)$ *of* $C + \hat{D}(\theta) + \theta \Delta_{\mathbf{v}}$ *equals 0 for a unique* $\theta = \alpha > 0$*, and we have, for* $i, j \in S$*,*

(5.9)
$$\lim_{x \to \infty} e^{\alpha x} P(M(\tau_x^+) = j | M(0) = i)$$
$$= \frac{1}{\eta^{(\alpha)}} \Big[ \mathbf{h}^{(\alpha)} \boldsymbol{\mu}^{(\alpha)} \Big( \Gamma(\alpha) - \frac{1}{\alpha}(C + D) \Big) \Big]_{ij},$$

*where* $\eta^{(\alpha)} = \boldsymbol{\mu}^{(\alpha)}(\hat{D}_{(1)}(\alpha) + \Delta_{\mathbf{v}})\mathbf{h}^{(\alpha)}$ *for the first derivative* $\hat{D}_{(1)}(\alpha)$ *of* $\hat{D}(\theta)$ *at* $\theta = \alpha$*, and* $\Gamma(\alpha)$ *is a* $S \times S$*-matrix such that*

$$(\Gamma(\alpha)^{+-}, \Gamma(\alpha)^{++}) = (0^{-+}, 0^{++}),$$
$$(\Gamma(\alpha)^{--}, \Gamma(\alpha)^{-+}) = -\mathbf{k}^{-}\boldsymbol{\pi}(\Delta_{\mathbf{v}} + \hat{D}'(0))$$
$$- \frac{1}{\alpha}(\alpha I^{--} - \mathbf{k}^{-}\boldsymbol{\pi}^{-})(\mathbf{k}^{-}\boldsymbol{\pi}^{-} - K^{(-)})^{-1}$$
$$\times (I^{--}, L^{-+})(C + D).$$

*In particular, using the notation* $\eta^{(0)} = E(Y(1)) < 0$*, for* $i \in S$*,*

(5.10) $$\lim_{x \to \infty} e^{\alpha x} P(\tau_x^+ < \infty | M(0) = i) = \frac{-\eta^{(0)} h_i^{(\alpha)}}{\eta^{(\alpha)}} (\boldsymbol{\mu}^{(\alpha)})^{-}\mathbf{k}^{-}.$$

PROOF. We first compute the denominator of (5.3). From (5.4), it is easy to see that

$$\boldsymbol{\nu}^{(\alpha)} \frac{d}{d\theta} \hat{H}(\theta) \Big|_{\theta=\alpha} \mathbf{h}^{(\alpha)} = \boldsymbol{\mu}^{(\alpha)}(\hat{D}_{(1)}(\alpha) + \Delta_{\mathbf{v}})\mathbf{h}^{(\alpha)},$$

since $(C + \hat{D}(\alpha) + \alpha \Delta_{\mathbf{v}})\mathbf{h}^{(\alpha)} = \mathbf{0}$. Thus the denominator is $\eta^{(\alpha)}$, which is clearly finite. This also implies condition (5.d). We next compute $\boldsymbol{\nu}^{(\alpha)} \int_0^{\infty} e^{\alpha u} \overline{G}(u)\, du$.

Then, (5.6) and (5.7) yield

$$(\mathbf{0}^-, (\boldsymbol{\nu}^{(\alpha)})^+) \int_0^\infty e^{\alpha u} \overline{G}(u)\, du$$
$$= \frac{1}{\alpha}(\mathbf{0}^-, (\boldsymbol{\mu}^{(\alpha)})^+)(C + \hat{D}(\alpha) + \alpha \Delta_{\mathbf{v}} - (C + D)).$$

Similarly, (5.6) and (5.8) yield

$$((\boldsymbol{\nu}^{(\alpha)})^-, \mathbf{0}^+) \int_0^\infty e^{\alpha u} \overline{G}(u)\, du$$
$$= \frac{1}{\alpha}((\boldsymbol{\mu}^{(\alpha)})^-, \mathbf{0}^+)(C + \hat{D}(\alpha) + \alpha \Delta_{\mathbf{v}} - (C + D))$$
$$- (\boldsymbol{\mu}^{(\alpha)})^- \mathbf{k}^- \boldsymbol{\pi}^- (I^{--}, L^{-+})(\Delta_{\mathbf{v}} + \hat{D}'(0))$$
$$- \frac{1}{\alpha}(\boldsymbol{\mu}^{(\alpha)})^- (\alpha I^{--} - \mathbf{k}^- \boldsymbol{\pi}^-)(\mathbf{k}^- \boldsymbol{\pi}^- - K^{(-)})^{-1}(I^{--}, L^{-+})(C + D).$$

Thus, (5.3) yields (5.9) since $\boldsymbol{\pi}^-(I^{--}, L^{-+}) = \boldsymbol{\pi}$ and $\boldsymbol{\mu}^{(\alpha)}(C + \hat{D}(\alpha) + \alpha \Delta_{\mathbf{v}}) = \mathbf{0}$. Finally, summing up (5.9) over all $j \in S$ yields (5.10), since $(C + D)\mathbf{e} = \mathbf{0}$ and

$$\boldsymbol{\pi}(\Delta_{\mathbf{v}} + \hat{D}_{(1)}(0))\mathbf{e} = E(Y(1)). \qquad \square$$

REMARK 5.2. When $S^+ = \varnothing$ and $v \equiv -1$, the results in Theorem 5.1 are fully compatible with those obtained in Theorem 3.1 of Miyazawa (2002).

It may be interesting to consider asymptotics of the probability that $Y(t)$ overshoots level $x$ continuously. To this end, we define $B$ of (5.2) as

(5.11) $$B(x) = \Delta_{\mathbf{v}}^{-1} \begin{pmatrix} -L^{-+} \\ I^{++} \end{pmatrix} \Delta_{\mathbf{v}}^{++} e^{-x \Delta_{\mathbf{c}}^{++} (\Delta_{\mathbf{v}}^{++})^{-1}}.$$

Then, using (5.6), we have

$$\boldsymbol{\nu}^{(\alpha)} \int_0^\infty e^{\alpha x} B(x)\, dx = \boldsymbol{\mu}^{(\alpha)} \begin{pmatrix} -L^{-+} \\ I^{++} \end{pmatrix} \Delta_{\mathbf{v}}^{++}.$$

Hence, we obtain the following result, which extends Corollary 4.9 of Asmussen (1994).

THEOREM 5.2. *Under the same conditions as in Theorem 5.1, for $i \in S$ and $j \in S^+$,*

(5.12)
$$\lim_{x \to \infty} e^{\alpha x} P(M(\tau_x^+) = j, Y(\tau_x^+) = x | M(0) = i)$$
$$= \frac{h_i^{(\alpha)}}{\eta^{(\alpha)}} \left[ \boldsymbol{\mu}^{(\alpha)} \begin{pmatrix} -L^{-+} \\ I^{++} \end{pmatrix} \Delta_{\mathbf{v}}^{++} \right]_j.$$

**6. Applications.** We briefly discuss applications of our results to a risk process and a fluid queue with extra jump inputs. Define a process $Z_x(t)$ for each $x \geq 0$ as

$$Z_x(t) = x - Y(t), \qquad t \geq 0.$$

Then, $Z_x(t)$ is a risk reserve process starting with reserve level $x$, and $-v(i)$ is the premium rate under background state $i \in S$, while a claim whose size has distribution $F_{ij}^{\mathrm{D}}$ arrives when the background state changes from $i$ to $j$ with rate $D_{ij}$. For the risk process $Z_x(t)$, a primary interest is the ruin probabilities and their asymptotics. The latter are obtained by Theorems 5.1 and 5.2. Numerical values of the ruin probabilities may be obtained taking Fourier transform of (4.1), using Theorem 4.1 and Corollary 4.1, and applying numerical inversion technique.

We next consider a fluid queue. Suppose that $Y(t)$ describes the net flow of a fluid queue with extra jump inputs. Then, the buffer content $V(t)$ at time $t$ is given by

$$V(t) = \sup_{0 \leq u \leq t} (Y(t) - Y(u)),$$

if $V(0) = 0$. As is well known, under the stability condition (5.1), which is assumed from now on, $V(t)$ and $M(t)$ converge jointly in distribution as $t \to \infty$. Denote a pair of random variables having this joint distribution by $(V, M)$. Then, the well-known formulation of Loynes (1962) yields

(6.1) $$P(V > x, M = i) = \pi_i P\left(\sup_{u \geq 0} \tilde{Y}(u) > x | M(0) = i\right).$$

Hence, this stationary probability is the hitting probability at level $x$ of the additive process $\tilde{Y}(t)$. Thus, converting (5.10) to the dual process, Theorem 5.1 yields the following:

THEOREM 6.1. *Under the conditions of Theorem 5.1, let $\alpha$ be the same one in the theorem, $Q^{(-)}$ be the minimal rate matrix satisfying the following matrix equation for a substochastic matrix $R^{+-}$:*

(6.2) $$-\begin{pmatrix} I^{--} \\ R^{+-} \end{pmatrix} Q^{(-)} = \Delta_{\mathbf{v}}^{-1} \int_0^\infty (C(du) + D(du)) \begin{pmatrix} I^{--} \\ R^{+-} \end{pmatrix} e^{uQ^{(-)}},$$

*and let $\boldsymbol{\beta}^-$ be the stationary probability vector of $Q^{(-)}$. Then, we have, for $i \in S$,*

(6.3) $$\lim_{x \to \infty} e^{\alpha x} P(V > x, M = i) = \frac{-\eta^{(0)} \mu_i^{(\alpha)}}{\eta^{(\alpha)}} \boldsymbol{\beta}^- (\mathbf{h}^{(\alpha)})^-.$$

# APPENDIX

**A.1. Proof of Lemma 5.2.** We first consider the case that $i \in S^+$. From (4.4), we have

$$\hat{H}_{ij}(\theta) = \delta_{ij} + \frac{1}{c(i) - v(i)\theta}(C_{ij} + \hat{D}_{ij}(\theta) + \theta[\Delta_{\mathbf{v}}]_{ij}).$$

Thus, we get (5.4). For $i \in S^-$, we first compute the moment generating function of the first term of (4.4).

$$\frac{-1}{v(i)}\left[\int_0^\infty e^{\theta x}\, dx \int_x^\infty e^{(y-x)K^{(-)}}(I^{--}, L^{-+})D(dy)\right]_{ij}$$

$$= \frac{-1}{v(i)}\left[(\theta I^{--} - K^{(-)})^{-1}\int_0^\infty (e^{y\theta I^{--}} - e^{yK^{(-)}})(I^{--}, L^{-+})D(dy)\right]_{ij}.$$

Hence, using the following equation that is obtained from (3.11),

$$-\int_0^\infty e^{yK^{(-)}}(I^{--}, L^{-+})D(dy)$$
$$= (I^{--}, L^{-+})C + (K^{(-)} - \theta I^{--})(I^{--}, L^{-+})\Delta_{\mathbf{v}} + \theta(I^{--}, L^{-+})\Delta_{\mathbf{v}},$$

we have

$$\frac{-1}{v(i)}\left[\int_0^\infty e^{\theta x}\, dx \int_x^\infty e^{(y-x)K^{(-)}}(I^{--}, L^{-+})D(dy)\right]_{ij}$$

$$= \delta_{ij} + \frac{-1}{v(i)}[(0^{--}, L^{-+})(-\Delta_{\mathbf{v}}) + (\theta I^{--} - K^{(-)})^{-1}(I^{--}, L^{-+})$$

$$\times (C + \theta\Delta_{\mathbf{v}} + \hat{D}(\theta))]_{ij}.$$

We next compute the moment generating function of the second term of (4.4). Similarly to the case that $i \in S^+$, we have

$$\frac{-1}{v(i)}\int_0^\infty e^{\theta x}\, dx \left(\int_0^x L_{ik}^{-+}\sum_{k \in S^+}(1 - e^{-c(k)(x-y)/v(k)})\frac{v(k)}{c(k)}U_{kj}(dy)\right)$$

$$= \frac{-1}{v(i)}[(0^{--}, L^{-+})\Delta_{\mathbf{v}}(I + \Delta_{\mathbf{c}-\theta\mathbf{v}}^{-1}(C + \hat{D}(\theta) + \theta\Delta_{\mathbf{v}}))]_{ij}.$$

**A.2. Derivation of (5.8).** Similarly to (5.7), (4.6) yields

$$\int_0^\infty e^{\alpha x}\overline{G}_{ij}(x)\, dx$$

(A.4)
$$= \frac{-1}{v(i)}\left[\int_0^\infty e^{\alpha x}dx \int_x^\infty dw \int_w^\infty e^{(y-w)K^{(-)}}(I^{--}, L^{-+})D(dy)\right.$$

$$\left. + \frac{1}{\alpha}(0^{--}, L^{-+})\Delta_{\mathbf{v}}\Delta_{\mathbf{c}-\alpha\mathbf{v}}^{-1}(C + \hat{D}(\alpha) + \alpha\Delta_{\mathbf{v}} - (C + D))\right]_{ij}.$$

We compute the integral in the bracket of (A.4) in the following way. Using the fact that $K^{(-)} - \alpha I^{--}$ and $\mathbf{k}^-\boldsymbol{\pi}^- - K^{(-)}$ are nonsingular, where the latter is obtained by Corollary 3.1, we have

$$\int_0^\infty e^{\alpha x}\,dx \int_x^\infty dw \int_w^\infty e^{(y-w)K^{(-)}}(I^{--},L^{-+})D(dy)$$

$$= \int_0^\infty dw \int_0^w e^{\alpha x}\,dx \int_w^\infty e^{(y-w)K^{(-)}}(I^{--},L^{-+})D(dy)$$

$$= \frac{1}{\alpha}\int_0^\infty dw \int_w^\infty (e^{\alpha w} - 1)e^{(y-w)K^{(-)}}(I^{--},L^{-+})D(dy)$$

$$= \frac{1}{\alpha}(K^{(-)} - \alpha I^{--})^{-1}\left(\int_0^\infty e^{yK^{(-)}}(I^{--},L^{-+})D(dy)\right.$$

$$\left. - (I^{--},L^{-+})\hat{D}(\alpha)\right)$$

$$- \frac{1}{\alpha}(\mathbf{k}^-\boldsymbol{\pi}^- - K^{(-)})^{-1}\int_0^\infty (y\mathbf{k}^-\boldsymbol{\pi}^- - (e^{yK^{(-)}} - I^{--}))$$

$$\times (I^{--},L^{-+})D(dy).$$

Using the relation (3.11), the first integral in the above equation is computed as

$$(\alpha I^{--} - K^{(-)})(I^{--},L^{-+})\Delta_{\mathbf{v}} - (I^{--},L^{-+})(C + \hat{D}(\alpha) + \alpha \Delta_{\mathbf{v}}),$$

while the second integral is computed as

$$(I^{--},L^{-+})(C + D) - (\mathbf{k}^-\boldsymbol{\pi}^- - K^{(-)})(I^{--},L^{-+})\Delta_{\mathbf{v}}$$

$$+ \mathbf{k}^-\boldsymbol{\pi}^-(I^{--},L^{-+})\left(\Delta_{\mathbf{v}} + \int_0^\infty yD(dy)\right).$$

Hence, using the fact that $\mathbf{k}^- = (\mathbf{k}^-\boldsymbol{\pi}^- - K^{(-)})^{-1}\mathbf{k}^-$, we have

$$\int_0^\infty e^{\alpha x}\,dx \int_x^\infty dw \int_w^\infty e^{(y-w)K^{(-)}}(I^{--},L^{-+})D(dy)$$

$$= \frac{1}{\alpha}\Big((\alpha I^{--} - K^{(-)})^{-1}(I^{--},L^{-+})(C + \hat{D}(\alpha) + \alpha \Delta_{\mathbf{v}})$$

$$- \mathbf{k}^-\boldsymbol{\pi}^-(I^{--},L^{-+})\left(\Delta_{\mathbf{v}} + \int_0^\infty yD(dy)\right)$$

$$- (\mathbf{k}^-\boldsymbol{\pi}^- - K^{(-)})^{-1}(I^{--},L^{-+})(C + D)\Big).$$

Substituting this formula into (A.4), we get (5.8).

**Acknowledgment.** The author is grateful to Søren Asmussen for pointing out the relation between (2.3) and (2.8) and related references Asmussen and Højgaard (1996) and Asmussen and Klüppelberg (1996).

Department of Information Sciences  
Tokyo University of Science  
Yamazaki 2641  
Noda  
Chiba 278-8510  
Japan  
e-mail: miyazawa@is.noda.tus.ac.jp